\documentclass[1p]{elsarticle}
\usepackage{amsmath,amsthm,amsfonts,amssymb,tabularx}
\usepackage[latin1]{inputenc}
\usepackage{graphicx}
\usepackage{color}
\usepackage{nameref}
\usepackage{pgf,tikz}
\usepackage{pdfsync}
\usetikzlibrary{arrows}
\usepackage{pgf,tikz}

\newtheorem{thm}{Theorem}

\newtheorem{defn}{Definition}
\newtheorem{pro}{Proposition}

\newcommand\lapal{(-\Delta)^{\alpha}}
\newcommand{\RN}{\mathbb{R}^N}
\newcommand\HH{\mathcal{H}}
\newcommand\KK{\mathcal{K}}
\newcommand\into{\int_{0}^T\int_{\mathbb{R}^N}}
\newcommand\intr{\int_{\mathbb{R}^N}}

\makeatletter
\def\ps@pprintTitle{%
  \let\@oddfoot\@empty
}
\makeatother

\begin{document}

\begin{frontmatter}

\textbf{\begin{center}
Partial Differential Equations
\end{center}}

\title{\bf Finite and infinite speed of propagation for  \\ porous medium equations with fractional pressure }

\author{Diana Stan}
\ead{diana.stan@uam.es}
\author{F\'{e}lix del Teso}
\ead{felix.delteso@uam.es}
\author{Juan Luis V{\'a}zquez\corref{cor1}}
\ead{juanluis.vazquez@uam.es}
\address{Departamento de Matem\'{a}ticas, Universidad
Aut\'{o}noma de Madrid,\\ Campus de Cantoblanco, 28049 Madrid, Spain}

\cortext[cor1]{Corresponding author.}

\begin{abstract}
We study a porous medium equation with fractional potential pressure:
$$
\partial_t u= \nabla \cdot (u^{m-1} \nabla p), \quad  p=(-\Delta)^{-s}u,
$$
for $m>1$, $0<s<1$ and $u(x,t)\ge 0$.  To be specific, the problem is posed for $x\in \RN$, $N\geq 1$, and $t>0$. The initial data $u(x,0)$ is assumed to be a bounded function with compact support or fast decay at infinity.
We establish existence of a class of weak solutions for which we determine whether,  depending on the parameter $m$, the property of compact support is conserved in time or not, starting from the  result of finite propagation known for $m=2$. We find that when $m\in [1,2)$ the problem has infinite speed of propagation, while for $m\in [2,\infty)$ it has finite speed of propagation. Comparison is made with other nonlinear diffusion models where the results are widely different.

\medskip

\noindent \textbf{R\'{e}sum\'e}

\textbf{Vitesse de propagation finie et infinie pour des \'{e}quations du milieu poreux avec une pression fractionnaire.} Nous \'{e}tudions une \'{e}quation du milieu poreux avec une pression potentielle fractionnaire: $
\partial_t u = \nabla \cdot (u^{m-1} \nabla p), \, p = (- \Delta)^{-s} u,
$
pour $m>1$, $0 <s <1 $ et $u (x,t) \ge 0 $. Le probl\`{e}me se pose pour $x \in \mathbb{R}^N$, $N \geq 1$ et $ t> 0 $. La donn\'{e}e initiale est suppos\'{e}e  born\'{e}e  avec support compact ou d\'{e}croissance rapide \`{a} l'infini. Lorsque le param\`{e}tre $ m $ est variable, on obtient deux comportements diff\'{e}rents comme suit: si  $m  \in [1,2) $ le probl\`{e}me a une  vitesse de propagation infinie, alors que pour  $m \in [2,\infty) $, elle a une vitesse  de propagation finie. On compare le r\'esultat avec le comportement d'autres mod\`eles de diffusion nonlin\'eaire qui est tr\`es diff\'erent.

\end{abstract}


\end{frontmatter}

\noindent \textbf{Version fran\c{c}aise abr\'{e}g\'e{e}}

Nous \'etudions un mod\`ele de diffusion non lin\'eaire avec pression nonlocale donn\'ee par
\begin{equation}\label{model1}
\partial_t u= \nabla \cdot (u^{m-1} \nabla p), \quad  p=\mathcal{K}(u),
\end{equation}
pour $m>1$ et $u(x,t)\geq0$. Le probl\`eme est pos\'e pour $x\in \RN$, $N\geq 1$, et $t>0$, et on se donne des conditions initiales
\begin{equation}\label{initialcondition}
u(x,0)=u_0(x), \quad x \in \RN,
\end{equation}
 o\`u $u_0:\RN \to [0,\infty)$ est born\'ee et \`a support compact ou \`a d\'ecroissance rapide \`a l'infini. La pression $p$ est li\'ee \`a $u$ par un op\'erateur potentiel lin\'eaire fractionnaire $p=\mathcal{K}(u)$, plus pr\'ecis\'ement   $\mathcal{K}=(-\Delta)^{-s} $  \`a $0<s<1 $, qui s'exprime par le noyau $K(x,y)=c|x-y|^{-(N-2s)}$ (c.\`a.d. un op\'erateur de Riesz).

Notre travail est motiv\'e par deux travaux r\'ecents qui consid\`erent l'\'equation \eqref{model1}  pour $m=2$. D'une part,  Caffarelli et V\'{a}zquez \'etudient dans  \cite{CaffVaz} le mod\`ele de diffusion nonlin\'eaire de type  milieu poreux avec des effets de diffusion non locaux:
\[
  \partial_t u= \nabla \cdot (u \nabla p), \ \ \ p=(-\Delta)^{-s}u. \tag{CV}\label{ModelCaffVaz}
\]
L'etude est compl\'et\'e dans les articles  \cite{CaffVaz2,CaffSorVaz}. D'autre part, Head \cite{Head} a d\'ecrit la dynamique des dislocations dans les cristaux vues comme un continuum, et a propos\'e l'\'equation (1.1) avec $s = 1/2$, $m = 2$ si la dimension est $N = 1$. La densit\'e de dislocations est de $u = v_x$ , donc $v$ r\'esout le ``probl\`eme int\'egr\'e''
$$
v_t+|v_x|(-\partial_{xx})^{1-s}v=0.
$$
Le mod\`ele sous cette forme a \'et\'e r\'ecemment \'etudi\'e par Biler, Karch et Monneau dans \cite{Biler2} o\`u ils prouvent l'existence d'une unique solution de viscosit\'e. Ils trouvent aussi une solution auto-similaire explicite et ils d\'ecrivent le comportement asymptotique.

Notre note a pour but de mettre en clair la propri\'et\'e de propagation finie pour les solutions du probl\`eme \eqref{model1}-\eqref{initialcondition} en dependance du param\`etre $m$. On sait d\'ej\`a que cette propri\'et\'e est valable pour $m = 2$. On prouve ici que elle est encore valable pour $m \in (2, \infty)$ tandis que pour $1 < m < 2$ les solutions non n\'egatives de cette \'equation sont strictement positives pour tout $x \in \mathbb{R}$ et tout $t > 0$. Ce dernier r\'esultat est prouv\'e rigoureusement en dimension $N = 1$ et sous des conditions sur les donn\'ees initiales. Notons finalement que les cas $m=1$ and $m=2$ \'etaient connus.

Tous ces points seront d\'evelopp\'ees dans l'article \cite{StanTesoVazquez}.

\section{Introduction}

We study a nonlinear diffusion model with nonlocal pressure given by
\begin{equation*}
\partial_t u= \nabla \cdot (u^{m-1} \nabla p), \quad  p=\mathcal{K}(u),
\end{equation*}
for $m>1$ and $u(x,t)\ge 0$. The problem is posed for $x\in \RN$, $N\geq 1$, and $t>0$, and we give initial conditions
\begin{equation*}
u(x,0)=u_0(x), \quad x \in \RN,
\end{equation*}
where $u_0:\RN \to [0,\infty)$  is bounded with compact support or fast decay at infinity.

The pressure $p$ is related to $u$ through a linear fractional potential operator  $p=\mathcal{K}(u)$. To be specific  $\mathcal{K}=(-\Delta)^{-s} $  for $0<s<1 $ with kernel $K(x,y)=c|x-y|^{-(N-2s)}$ (i.e. a Riesz operator). Moreover, $\mathcal K$ is a self-adjoint operator with $\mathcal K=\mathcal H^2$.

Our work is motivated by two recent works that consider equation \eqref{model1} for $m=2$. On one hand, Caffarelli and V\'{a}zquez have studied \eqref{model1} for $m=2$ in \cite{CaffVaz} as a nonlinear diffusion equation of porous medium type with nonlocal diffusion effects
\[
  \partial_t u= \nabla \cdot (u \nabla p), \ \ \ p=(-\Delta)^{-s}u.
\]
The study of this model has been pursued \cite{CaffVaz2,CaffSorVaz}. The properties of boundedness of the solutions, $C^\alpha$ regularity, existence of self-similar solutions and asymptotic  behavior are established. The latter paper points out that the application to problems in particle systems with long range interactions \cite{GLP} leads to more general equations of the form $u_t=\nabla\cdot(f(u)\nabla {\cal K} g(u))$ with convenient monotone functions $f$ and $g$, thus motivating the interest in our present model.

On the other hand, a similar model was proposed by Head \cite{Head}  to describe the dynamics of dislocations in crystals seen as a continuum. When the space dimension is $N=1$, the model becomes equation \eqref{model1} $s=1/2$ and $m=2$. The dislocation density is $u=v_x$, therefore $v$ solves the ``integrated problem''
$$
v_t+|v_x|(-\partial_{xx})^{1-s}v=0.
$$
The model in this form has been recently studied by Biler, Karch and Monneau in \cite{Biler2} where they prove existence of a  unique viscosity solution.


\section{Main results}

We first propose a definition of solution and establish the existence and main properties of the solutions.

\begin{defn}
We say that $u$ is a weak solution of $\eqref{model1} - \eqref{initialcondition}$ in $Q_T=\RN \times (0,T)$ with nonnegative initial data $u_0\in L^1(\RN)$ if \ (i)  $u\in L^1(Q_T)$,
 (ii) $\KK(u) \in L^1(0,T: W^{1,1}_{loc}(\RN))$, \
(iii) $u^{m-1}\nabla \KK(u) \in L^1(Q_T)$, and (iv)
\begin{equation}\label{model1weak}
\into u \phi_tdxdt -\into u^{m-1}\nabla \KK (u) \nabla \phi dxdt+ \intr u_0(x)\phi(x,0)dx=0
\end{equation}
holds for every test function $\phi$ in $Q_T$ such that $\nabla \phi$ is continuous, $\phi$ has compact support in $\RN$ for all $t\in (0,T)$ and vanish near $t=T$.
\end{defn}

Let state now the existence results as well as the basic properties for the solutions of problem \eqref{model1}.

\begin{thm}\label{Thm1}
Let  $m\in (1,2)$ (and $s\in (0,1/2)$ if $N=1$). Let $u_0\in L^1(\RN)\cap L^\infty (\RN)$.
Then there exists a weak solution $u$ of equation \eqref{model1} with initial data $u_0$.  Moreover, $u$ has the following properties:
\begin{enumerate}
\item  \textbf{(Regularity)} $u\in C([0,\infty): L^1(\RN))$, $u \in L^\infty (\RN\times (0,T))$, $\nabla\mathcal{H}(u) \in L^2(Q)$ .
\item \textbf{(Conservation of mass)} For all $t>0$ we have
$\displaystyle{
\int_{\RN}u(x,t)dx=\int_{\RN}u_0(x)dx.
}$
\item \textbf{($L^{\infty}$ estimate) } For all $t>0$ we have  $||u(\cdot,t)||_\infty\leq ||u_0||_\infty$.

\item \textbf{(Energy estimate)} For all $t>0$ the following estimate holds
$$\displaystyle{
C\int_0^t  \int_{\RN} |\nabla \HH(u)|^2dxdt  +\int_{\RN}u(t)^{3-m}dx =\int_{\RN}u_0^{3-m}dx}, \qquad C=(2-m)(3-m)>0\,.
$$
\end{enumerate}
\end{thm}

\begin{thm}\label{Thm2}
Let $m\in [2,\infty)$. Let $u_0\in L^1(\RN)\cap L^\infty (\RN)$ be such that
\begin{equation}
0\le u_0(x)\le Ae^{-a|x|} \text{ for some }A,a>0.
\end{equation}
Then there exists a weak solution $u$ of equation \eqref{model1} with initial data $u_0$
such that $u$ satisfies the properties $1,2,3$ of Theorem \ref{Thm1}.  Moreover, the solution
decays exponentially in $|x|$ and enjoys  suitable energy estimates, somewhat different from the ones of Theorem \ref{Thm1}.
\end{thm}

For $m=2$ see \cite{CaffVaz}.  The following is our most important contribution, which deals with the property of finite propagation of the solutions just constructed depending on the value of $m$.

\begin{thm}
Let $m\in (2,\infty)$. The solution to problem \eqref{model1}-\eqref{initialcondition} has the property of \textbf{finite speed of propagation} in the sense that, if  $u_0$ is compactly supported, then for any $t>0$, $u(\cdot,t)$ is also compactly supported. On the other hand, if $N=1$, $m\in (1,2)$ and $s\in (0,1/2)$, the solution $u$ has \textbf{infinite speed of propagation} in the sense that: even if the initial data is compactly supported, for any $t>0$ and any $R>0$, the set $\mathcal{M}_{R,t}=\{x: |x|\ge R,\  u(x,t)>0\}$ has positive measure.
\end{thm}

\noindent\textbf{Remarks.} (i) Finite propagation for the case $m=2$ has been proved in \cite{CaffVaz}.

\noindent (ii) A main difficulty in the work to be done is the possible lack of uniqueness and comparison of the solutions, already noticed in \cite{CaffVaz}. This is reflected in the indirect statement of the last part of the theorem.

\noindent (iii) Since $\nabla \cdot (u^{m-1} \nabla (-\Delta)^{-s}u)$ tends to $-(-\Delta)^{1-s}u $ as $m\to 1$,  equation \eqref{model1} becomes $u_t+(-\Delta)^{1-s}u=0$, known as the fractional heat equation, which has infinite speed of propagation. This propagation property is inherited by more general diffusion models as $u_t+(-\Delta)^s u^m=0$, called the fractional porous medium equation, which has been  studied in \cite{PQRV1,PQRV2,VazquezBarenblattFractPME}. Therefore, a change in the behavior of the solutions for  some $m> 1$ was expected. We also motivate the result by numerical computations based on the scheme proposed by Teso and V\'azquez in \cite{Teso2013,TesoVaz2013}. We give in the graph below (Figure \ref{fig1}) a description of the result on finite propagation for different related models of nonlinear diffusion with or without fractional effects, see \cite{VazSurveyFractional, PQRV2, Biler2}.
\definecolor{qqqqff}{rgb}{0,0,1}
\definecolor{ffqqtt}{rgb}{1,0,0.2}
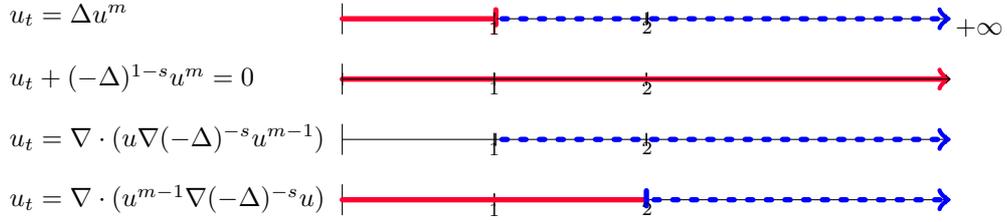
\begin{figure}
\centering
\begin{tikzpicture}[line cap=round,line join=round,x=1cm,y=0.4cm]
\clip(-2.88,0.58) rectangle (11.04,8.69);
\draw (2,8.5)-- (2,7.5);
\draw [line width=2pt,color=ffqqtt] (4.02,8.3)-- (4.02,7.8);
\draw (6,8.2)-- (6,7.8);
\draw (4,8.2)-- (4,7.8);
\draw (6,8.2)-- (6,7.8);
\draw [->] (2,8) -- (10,8);
\draw (2,6.5)-- (2,5.5);
\draw (6,6.2)-- (6,5.8);
\draw (4,6.2)-- (4,5.8);
\draw (6,6.2)-- (6,5.8);
\draw [->,line width=2pt,color=ffqqtt] (2,6) -- (10,6);
\draw (2,4.5)-- (2,3.5);
\draw (4,4.2)-- (4,3.8);
\draw (6,4.2)-- (6,3.8);
\draw (4,4.2)-- (4,3.8);
\draw (6,4.2)-- (6,3.8);
\draw [->] (2,4) -- (10,4);
\draw (2,2.5)-- (2,1.5);
\draw (4,2.2)-- (4,1.8);
\draw (6,2.2)-- (6,1.8);
\draw (4,2.2)-- (4,1.8);
\draw (6,2.2)-- (6,1.8);
\draw [->] (2,2) -- (10,2);
\draw [line width=2pt,color=ffqqtt] (2,8)-- (4,8);
\draw [line width=2pt,color=ffqqtt] (2,2)-- (6,2); 
\draw (4,8.2)-- (4,7.8);
\draw (4,6.2)-- (4,5.8);
\draw (2,6.2)-- (2,5.8);
\draw (4,6.2)-- (4,5.8);
\draw (2,6.2)-- (2,5.8);
\draw [line width=2pt,color=qqqqff] (6,2.3)-- (6,1.8); 
\draw (4.02,4.2)-- (4.02,3.8);
\draw [->] (2,6) -- (10,6);
\draw (9.94,8.3) node[anchor=north west] {$+\infty$};
\draw [->,line width=2pt,dash pattern=on 2pt off 5pt,color=qqqqff] (6.1,2) -- (10,2);
\draw [->,line width=2pt,dash pattern=on 2pt off 5pt,color=qqqqff] (4.1,8) -- (10,8);
\draw [->,line width=2pt,dash pattern=on 2pt off 5pt,color=qqqqff] (4.1,4) -- (10,4);
\draw (6.02,8.22)-- (6.02,7.82);
\draw (-2.5,8.8) node[anchor=north west] {$u_t=\Delta u^m$};
\draw (-2.5,6.8) node[anchor=north west] {$u_t+(-\Delta)^{1-s}u^m=0$};
\draw (-2.5,4.8) node[anchor=north west] {$u_t=\nabla \cdot (u\nabla (-\Delta)^{-s}u^{m-1})$};
\draw (-2.5,2.8) node[anchor=north west] {$u_t=\nabla \cdot (u^{m-1}\nabla (-\Delta)^{-s}u)$};
\begin{scriptsize}
\draw[color=ffqqtt] (4,7.67) node {1};
\draw[color=black] (6.01,7.71) node {2};
\draw[color=black] (6.01,5.7) node {2};
\draw[color=black] (4,3.66) node {1};
\draw[color=black] (6.01,3.7) node {2};
\draw[color=black] (4,1.66) node {1};
\draw[color=black] (6.01,1.69) node {2};
\draw[color=black] (4,7.67) node {1};
\draw[color=black] (4,5.7) node {1};
\end{scriptsize}
\end{tikzpicture}
\caption{Intervalles de l'exposant $m\in (0,\infty).$ La ligne rouge signifie propagation infinie, la ligne pointill\'ee bleue est la propagation finie. Ranges of exponent $m\in (0,\infty).$ Red line means infinite propagation, blue dotted line is finite propagation. } \label{fig1}
\end{figure}


\section{Existence of weak solutions}

We make an approach to problem \eqref{model1} based on regularization, elimination of the degeneracy and reduction of the spatial domain.  Once we have solved the approximate problems, we derive estimates that allow us to pass to the limit in all the steps one by one, to finally obtain a weak solution of the original problem. In doing this we follow the outline of (\cite{CaffVaz}). Specifically, for small $\epsilon, \delta, \mu \in (0,1)$ and $R>0$ we consider the following initial boundary value problem
posed in $Q_{T,R}=\{x\in B_R(0), \ t\in (0,T)\}$
\begin{equation}\label{model1Aprox}
\left\{
\begin{array}{ll}
u_t= \delta \Delta u +\nabla ((u+\mu)^{m-1} \nabla \KK_\epsilon (u))&\text{for } (x,t)\in Q_{T,R}\\
u(x,0)=\widehat{u}_0(x) &\text{for } x \in B_R(0)\\
u(x,t)=0 &\text{for } x\in \partial B_R(0), \ t\geq 0,
\end{array}
\right.
\end{equation}
where $\KK_\epsilon$ is a smooth approximation of $(-\Delta)^{-s}$. Also, $\mathcal{K}_\epsilon$ is a self-adjoint operator and we denote by $\mathcal{K}_\epsilon=\mathcal{H}_\epsilon^2$. The basic tool of the proof of existence is the following energy estimate for the approximated problem \eqref{model1Aprox}:
\begin{equation}\label{energy1}
\int_{B_R} F_\mu(u(t))dx+\delta\int_0^t \int_{B_R} \frac{|\nabla u|^2}{d_\mu(u)}dxdt+\int_0^t  \int_{B_R} |\nabla   \HH_\epsilon(u)|^2dxdt=\int_{B_R}F_\mu(\widehat{u}_0)dx.
\end{equation}
for all $0<t<T$.  The explicit formula for $F_\mu$  is:
\begin{equation*}\displaystyle{
F_\mu(u)=
\left\{
\begin{array}{ll}
 \frac{1}{(2-m)(3-m)}[(u+\mu)^{3-m}-\mu^{3-m}]-\frac{1}{2-m}\mu^{2-m}u,&\text{for } m\not=2,3\\[3mm]
-\log\left(1+(u/\mu  \right)) + u/\mu ,&\text{for }m=3\\[2mm]
(u+\mu) \log\left(1+(u/\mu)\right)-u, &\text{for } m=2.
\end{array}
\right.}
\end{equation*}
For $m=2$ see \cite{CaffVaz}.

\subsection{Uniform Bounds in the case $m\in (1,2)$}
We obtain uniform bounds in all parameters $\epsilon, R, \delta,\mu$ for the energy estimate \eqref{energy1}, that allows us to pass to the limit and obtain a solution of the original problem \eqref{model1}. More exactly,
$$
\frac{1}{(m-2)(3-m)}\int_{\RN}\left[(u_0+\mu)^{3-m}-\mu^{3-m}\right]dx \leq \frac{(||u_0||_\infty+1)^{2-m}}{m-2}\int_{\RN} u_0dx.
$$

\subsection{Exponential tail in the case $m\in (2,\infty)$}
This case is more delicate since  the term $\int_{B_R} F_\mu(u(t))dx$ can not be easily uniformly controlled in $\mu>0$.

In \cite{CaffVaz}, when $m=2$, the authors prove an exponential tail control of the approximate solution by using a comparison method with a suitable family of barrier functions, called true supersolutions.  Their proof can be adapted to the case $m\in (2,\infty)$ with a series of technical modifications caused by the power $u^{m-1}.$

We observe that for $m\in (2,3)$
\[
\frac{1}{(2-m)(3-m)}\int_{\RN}\left[(u(t)+\mu)^{3-m}-\mu^{3-m}\right]dx \leq \frac{1}{(2-m)(3-m)}\int_{\RN} u(t)^{3-m}dx,
\]
which is finite due to the tail control result.

\section{Finite propagation for $m\in (2,\infty)$}

 Assume $u$ is a bounded solution, $0\le u\le L$, of Problem \eqref{model1} as constructed before.  Assume that $u_0$ has compact support. Then we prove that $u(\cdot,t)$ is compactly supported for all $t>0$. More precisely, if $u_0$ is below the ''parabola-like'' function
$$U_0(x)=a(|x|-b)^2,$$
for some $a,b>0$, with support in the ball $B_b(0)$, then for $C$ is large enough, we prove that
$$
u(x,t) \le a(Ct-(|x|-b))^2.
$$
Since the problem does not satisfy the comparison principle, the proof of the above result is made by arguing at the first point in space and time  where $u(x,t)$ touches the parabola $U$ from below. A contradiction is obtained as in the typical viscosity method if $a$ and $b$ suitable chosen, even if there is generally no comparison principle.

\section{Infinite propagation for $m\in (1,2)$}

The properties of the Model \eqref{model1} in dimension $N=1$ can be obtained via a study of the properties of the integrated problem.
We define $v(x,t)=\int_{-\infty}^x u(y,t)dy$. Then $v$ is a solution to the Problem
\begin{equation}\label{IntegEq}
  \left\{ \begin{array}{ll}
  \partial_t v= -|v_x|^{m-1}\lapal v &\text{for }x \in \mathbb{R}\text{ and } t>0 , \\[2mm]
  v(x,0)=v_0(x)&\text{for } x \in \mathbb{R},
    \end{array}
    \right.
\end{equation}
where the initial data is given by $v_0(x)=\int_{-\infty}^x u_0(y)dy $. The exponents $\alpha$ and $s$ are related by
\begin{equation}\label{alphas}
\alpha=1-s.
\end{equation}
The sketch of the proof of the infinite speed of propagation is as follows.

We define the notions of viscosity sub-solution, super-solution and solution in the sense of Crandall-Lions \cite{Crandall}. The definition will be adapted to our problem by considering the time dependency and also the nonlocal character of the Fractional Laplacian operator.

The existence of a unique viscosity solution to Problem \eqref{IntegEq} follows as in \cite{Biler2}. The standard comparison principle for viscosity solutions also holds true (see \cite{ImbertHomog, Jakobsen}).

We give now our extended version of parabolic comparison principle, which represents an important instrument when using barrier methods. A similar parabolic comparison has been proved in \cite{CabreRoque} and has been used for instance in \cite{CabreRoque,KPPStanVazquez}. Let $\Omega \subset \mathbb{R}$ possibly unbounded.

\begin{pro}\label{ComparPrinc1}
Let $m>1$, $\alpha \in (0,1)$.  Let $v$ be a viscosity solution of Problem \eqref{IntegEq}. Let  $\Phi: \mathbb{R}\times[0,\infty)\to \mathbb{R}$ such that $\Phi \in C^2(\Omega\times (0,T))$. Assume that
\begin{itemize}
  \item $\Phi_t +|\Phi_x|^{m-1}\lapal\Phi< 0 $ for $x\in \Omega$, $t\in [0,T]$.
  \item $\Phi(x,0) < v(x,0)$  for all $x\in \mathbb{R}^N$ (comparison at initial time).
  \item $\Phi(x,t) < v(x,t)$  for all $x \in \RN \setminus \Omega$ and $t\in (0,T)$ (comparison on the parabolic boundary).
\end{itemize}
Then $\Phi(x,t) \le v(x,t)$ for all $x \in \mathbb{R}$, $t\in (0,T).$
\end{pro}

We construct a lower barrier of the form
\begin{equation}\label{barrier1}
\Phi_\epsilon(x,t)=(t+\tau)^{b \gamma}\left((|x|+\xi)^{-\gamma}+G(x)\right)-\epsilon, \quad t\geq0, \ x\in \mathbb{R}.
\end{equation}
where $\gamma=(2\alpha+m)/(2-m)$ is the exponent of a kind of self-similar solution of this equation. The function $G$ is such that, given any two constants $C_1>0$ and $C_2>0$, we have that\\
 $\bullet$ (G1) $G$ is compactly supported in the interval $(-x_0,\infty)$.\\
 $\bullet$ (G2) $G(x) \le C_1$ for all $x\in \mathbb{R}$.\\
 $\bullet$ (G3) $(-\Delta)^s G(x) \le -C_2|x|^{-(1+2s)}$ for all $x<x_0$.

We also prove the existence of this kind of function. To finish, we show that, and for any $t_1>0$ and any $x_1\in \mathbb{R}$ we have that $\Phi_\epsilon$ satisfies the hypothesis of Proposition  \ref{ComparPrinc1} and  $\Phi_\epsilon(x_1,t_1)>0$ for a suitable choice of the parameters $\epsilon, \xi$ and $\tau$. A similar lower barrier has been constructed in \cite{KPPStanVazquez}.

All these results will be developed in \cite{StanTesoVazquez}.

\section{Comments and Open Problems}

\noindent $\bullet$ \textbf{Explicit solutions.} Y. Huang reports  \cite{Huang2013} the explicit expression of the Barenblatt solution for the special value of  $m$, $m_{ex} = (N+6s-2)/(N+2s)$. The profile is given by
\begin{equation*}
F_M(y) =\lambda\,(R^2 + |y|^2)^{-(N+2s)/2}
\end{equation*}
where the two constants $\lambda$ and $R$ are determined by the total mass $M$ of the solution  and the parameter $\beta$. Note that for $s=1/2$ we have $m_{ex}=1$, and the solution corresponds to the linear case, $u_t=(-\Delta)^{1/2} u$, $F_{1/2}(r)=C(a^2+r^2)^{-(N+1)/2}$.

\noindent Different generalizations of model \eqref{ModelCaffVaz} are worth studying:

\noindent$\bullet$ Changing-sign solutions for the problem $\displaystyle{\partial_t u= \nabla \cdot (|u| \nabla p), \quad  p=(-\Delta)^{-s}u.}$

\noindent $\bullet$ Starting from the Problem \eqref{ModelCaffVaz}, an alternative is to consider the problem
\begin{equation*}\label{modelKarch}
u_t=\nabla.(|u|\nabla (-\Delta)^{-s}(|u|^{m-2}u)), \quad x \in \RN, \ t>0,
\end{equation*}
with $m>1$. This problem has been studied by Biler, Imbert and Karch in \cite{BilerCRAS, Biler1}. They construct explicit compactly supported self-similar
solutions which generalize the Barenblatt profiles of the PME. They do not prove the finite propagation of a general solution.

\noindent $\bullet$ We should consider combining the above models into $ \displaystyle \partial_tu=\nabla(|u|^{m-1}\nabla p),  \quad p=(-\Delta)^{-s}u. $ \\ When $s=0$ and $m=2$ we obtain the signed porous medium equation  $ \partial_tu=\Delta( |u|^{m-1}u)$.

{\sc Acknowledgment.} Work partially supported by Spanish Project MTM2011-24696 (Spain). The second author is also supported by a FPU grant from MECD, Spain.


\end{document}